\documentclass[a4paper,12pt]{article}
\usepackage{amsmath,amscd,amssymb,french}
\newcommand{\qed}{\hfill$\square$}

\begin{document}
\begin{center}
\large{\textbf{Structures de Contact sur les\\ 
Vari\'et\'es Alg\'ebriques de dimension 5}}
\end{center}
\vspace{0.5cm}
\centerline{St\'ephane DRUEL}
\begin{center}
DMI-\'Ecole Normale Sup\'erieure\\
45 rue d'Ulm\\
75005 PARIS\\
e-mail: \texttt{druel@clipper.ens.fr}
\end{center}
\vspace{1cm} 
\textbf{Introduction}\\
\newline
\indent Une \textit{structure de contact} sur une vari\'et\'e
alg\'ebrique lisse est la donn\'ee d'un sous-fibr\'e
$D\subset\mathcal{T}_{X}$ de rang dim$X-1$ de sorte que la forme
$\mathcal{O}_{X}$-bilin\'eaire sur $D$ \`a valeurs dans le fibr\'e en droites
$L=\mathcal{T}_{X}/D$ d\'eduite du crochet de Lie sur
$\mathcal{T}_{X}$ soit non d\'eg\'en\'er\'ee en tout point de
$X$. Cela entra\^{\i}ne que $X$ est de dimension impaire $2n+1$ et que le
fibr\'e canonique $K_{X}$ est isomorphe \`a
$L^{-1-n}$. On peut aussi d\'efinir la structure de contact par la
donn\'ee d'un \'el\'ement $\theta\in
H^{0}(X,\Omega_{X}^{1}\otimes L)$, la \textit{forme de contact},
tel que $\theta\wedge(d\theta)^{n}$ soit partout non nul.\\
\indent Soit $\mathfrak{g}$ une alg\`ebre de Lie simple. Son 
groupe adjoint $G$ agit sur $\mathbb{P}(\mathfrak{g})$ et n'a qu'une
seule orbite ferm\'ee; on montre que celle-ci admet une structure de
contact ([4] prop.2.6). Ce sont des vari\'et\'es de Fano homog\`enes
dont le groupe des automorphismes de contact a
pour alg\`ebre de Lie $\mathfrak{g}$. On parlera de la vari\'et\'e
de contact homog\`ene de type $\mathfrak{g}$. Les fibr\'es
$\mathbb{P}_{Y}(\mathcal{T}_{Y})$, o\`u $Y$ est une vari\'et\'e lisse,
fournissent d'autres exemples de vari\'et\'es de contact. En dimension 3,
Y-G.Ye a montr\'e que ce sont les seules ([17]). D'autres auteurs \'etudient
les structures de contact sur les vari\'et\'es de Fano ([4], [11]), mais
les r\'esultats ne sont ici encore que partiels.\\
\indent Le r\'esultat principal de ce travail est le \\
\newline
\textbf{Th\'eor\`eme}.$-$\textit{Soit $X$ est une vari\'et\'e projective
lisse de dimension 5 munie d'une structure de contact. Alors $X$ est
l'une des vari\'et\'es pr\'ec\'edentes sauf si le fibr\'e canonique
$K_{X}$ est num\'eriquement effectif et $\kappa(X)=-\infty$.}\\
\newline
\indent Notons que le dernier cas ne devrait pas se produire  par la conjecture
d'Abondance ([10]).\\ 
\newline
\textbf{Remerciements}.$-$Je tiens \`a exprimer toute ma gratitude
\`a A.Beauville pour m'avoir sugg\'er\'e ce probl\`eme et pour l'aide
qu'il m'a apport\'e.\\
\newline
\textbf{1. Rappels}\\
\newline
\indent Soit $X$ une vari\'et\'e projective lisse sur le corps
$\mathbb{C}$ des nombres complexes. Le produit d'intersection entre
1-cycles et diviseurs met en dualit\'e les deux espaces vectoriels
r\'eels :
$$N_{1}(X)=(\{\text{1-cycles}\}/\equiv)\otimes\mathbb{R}\text{ et } 
N^{1}(X)=(\{\text{diviseurs}\}/\equiv)\otimes\mathbb{R},$$
\noindent o\`u $\equiv$ d\'esigne l'\'equivalence num\'erique. La
dimension commune de ces espaces vectoriels est appel\'ee le 
\textit{nombre de
Picard} de $X$. On consid\`ere le c\^one $NE(X)\subset N_{1}(X)$
engendr\'e par les classes des 1-cycles effectifs. Une \textit{raie
extr\'emale} est une demi-droite $R$ dans $\overline{NE}(X)$, 
adh\'erence de $NE(X)$ dans $N_{1}(X)$, v\'erifiant $K_{X}.R^{*}<0$  et telle que pour tout $Z_{1},Z_{2}\in\overline{NE}(X)$, si
$Z_{1}+Z_{2}\in R$ alors $Z_{1},Z_{2}\in R$.
 Une \textit{courbe rationnelle extr\'emale} est une 
courbe rationnelle irr\'eductible $C$ telle que
$\mathbb{R}^{+}[C]$ soit une raie extr\'emale et
$-K_{X}.C\le\text{dim}X+1$. Le premier r\'esultat de la th\'eorie de Mori est que \textit{toute
raie extr\'emale est engendr\'ee par une courbe rationnelle
extr\'emale.} Le second r\'esultat fondamental est que 
\textit{toute raie 
extr\'emale $R$ admet une contraction}, c'est-\`a-dire qu'il existe une
vari\'et\'e projective normale $Y$ et un morphisme
$X\overset{\phi}{\longrightarrow}Y$, surjectif \`a fibres connexes,
contractant les courbes irr\'eductibles $C$ telles que $[C]\in R$
(th\'eor\`eme de Kawamata-Shokurov).\\
\indent Rappelons un r\'esultat fondamental de 
J.Wisniewski ([15], [16]) sur le 
lieu exceptionnel d'une contraction extr\'emale. Soit $F$ une
composante irr\'eductible d'une fibre non triviale d'une contraction
\'el\'ementaire associ\'ee \`a la raie extr\'emale $R$. Nous appelons
\textit{lieu de $R$}, le lieu des courbes dont la classe d'\'equivalence
num\'erique appartient \`a $R$. On a alors l'in\'egalit\'e :
$$\text{dim}F+\text{dim}(\text{lieu de }R)\ge\text{dim}X+\ell(R)-1,$$
\noindent o\`u $\ell(R)$ d\'esigne la \textit{longueur} de la raie extr\'emale
$R$ :
$$\ell(R)=\text{inf}\{-K_{X}.C_{0}| C_{0} \text{ \'etant une courbe
rationnelle et }C_{0}\in R\}.$$
\indent En particulier :
$$2\text{dim}(\text{lieu de }R)\ge\text{dim}X+\ell(R)-1.$$
\indent Nous terminons ces rappels par un th\'eor\`eme de
structure ([2], [3]). Consid\'erons une contraction extr\'emale
$X\overset{\phi}{\longrightarrow}Y$ d'une vari\'et\'e projective
lisse. Soit $L$ un fibr\'e inversible $\phi$-ample et $r\ge1$ un
entier. On dit que
$K_{X}+rL$ \textit{supporte la contraction} $\phi$ si ce fibr\'e est
trivial sur
les fibres de $\phi$. Soit $F=\phi^{-1}(y)$ une fibre de
$\phi$ munie de sa structure de sch\'ema r\'eduite. On suppose qu'il
existe un ouvert de $Y$ contenant $y$ tel que les fibres de $\phi$ au
dessus de cet ouvert soient toutes de
dimension au plus $\text{dim}(F)$ :
\begin{enumerate}
\item si $\text{dim}(F)\le r-1$ alors $Y$ est lisse en $y$ et $\phi$
est un fibr\'e projectif au voisinage de $F$,
\item si $\text{dim}(F)=r$ alors $Y$ est lisse au voisinage de
$y$ et :
\begin{enumerate}
\item si $\phi$ est birationnel alors $\phi$ est
l'\'eclatement d'une sous vari\'et\'e lisse de $Y$ de codimension
$r+1$,
\item si $\text{dim}(Y)=\text{dim}(X)-r$ alors $\phi$ est un fibr\'e en
quadriques,
\item si $\text{dim}(Y)=\text{dim}(X)-r+1$ alors $r\le \text{dim}(X)/2$ et
$F=\mathbb{P}^{r}$.
\end{enumerate}
\end{enumerate}
$\ $
\newline
\textbf{2. Preuve du th\'eor\`eme}\\
\newline
\textbf{Lemme}.$-$\textit{Soit $X$ une vari\'et\'e
projective lisse de dimension $2n+1$ munie d'une structure de contact
d\'efinie par la forme $\theta\in H^{0}(X,\Omega_{X}^{1}\otimes L)$. Soit
$Y\subset X$ une sous-vari\'et\'e analytique complexe lisse telle que 
la restriction de la forme de contact \`a $Y$ soit identiquement
nulle. Alors la dimension de $Y$ est au plus $n$.}\\
\newline
\textit{D\'emonstration}.$-$ On v\'erifie par un calcul en coordonn\'ees
locales que pour tout $y\in Y$, l'espace vectoriel
$\mathcal{T}_{Y}(y)\subset D(y)$ est un sous espace totalement
isotrope pour la forme altern\'ee de contact qui, par hypoth\`ese, est
non d\'eg\'en\'er\'ee.\qed\\
\newline
\textbf{Proposition 1}.$-$\textit{Soit $X$ une vari\'et\'e de Fano de
dimension 5 munie d'une structure de contact. On suppose que
$b_{2}(X)=1$. Alors $X$ est soit isomorphe \`a l'espace projectif
$\mathbb{P}^{5}$ soit \`a la vari\'et\'e de contact homog\`ene de type
$G_{2}$.}\\
\newline
\textit{D\'emonstration}.$-$Puisque $b_{2}(X)=1$, le groupe de Picard
de $X$ est un
$\mathbb{Z}$-module libre de rang 1. Rappelons que nous avons la
formule $K_{X}=-3L$. Il en r\'esulte que soit $L$ engendre
$\text{Pic}(X)$ soit $L=2L_{0}$ o\`u $L_{0}$ est un g\'en\'erateur du
groupe de Picard de $X$. Dans ce dernier cas,
$X$ est isomorphe \`a l'espace projectif complexe $\mathbb{P}^{2n+1}$
par le crit\`ere de Kobayashi-Ochiai ([9]).\\
\indent Il nous reste donc \`a traiter  le cas o\`u
$L$ est un g\'en\'erateur de $\text{Pic}(X)$. Dans ce cas, $X$ est une
vari\'et\'e de Mukai et le fibr\'e $L$ est tr\`es ample ([14] prop. 1,
[12]). En effet, lorsque $X$ est un rev\^etement double de
$\mathbb{P}^{5}$ ou d'une quadrique lisse de dimension $5$, on
v\'erifie que $H^{0}(X,\Omega^{1}_{X}\otimes L)=0$ et donc $X$ n'a
aucune structure de contact. Par suite,
$X$ est homog\`ene ([4] cor. 1.8) et isomorphe 
\`a la vari\'et\'e de contact homog\`ene de type $G_{2}$ ([5]).\qed\\
\newline
\textbf{Th\'eor\`eme 1}.$-$\textit{Soit $X$ une vari\'et\'e projective
lisse de dimension 5 munie d'une structure de contact. On suppose que
le fibr\'e canonique n'est pas num\'eriquement effectif. Alors $X$
est soit isomorphe \`a l'espace 
projectif $\mathbb{P}^{5}$, soit \`a
$\mathbb{P}_{Y}(\mathcal{T}_{Y})$ o\`u $Y$ est une vari\'et\'e lisse
de dimension 3, soit \`a la vari\'et\'e de contact homog\`ene de type
$G_{2}$.}\\
\newline
\textit{D\'emonstration}.$-$La preuve de ce th\'eor\`eme repose sur
l'\'etude des contractions ext\'emales de $X$. Soit $R$ une raie
extr\'emale de $X$. Puisqu'on a la formule $K_{X}=-3L$,
$\ell(R)=3\text{ ou }6$. Dans le dernier cas, $X$ est de Fano et
$b_{2}(X)=1$ ([15]) et la proposition 1 permet de conclure.\\
\indent Il nous reste \`a traiter le cas o\`u $\ell(R)=3$. Notons
$X\overset{\phi}{\longrightarrow}Y$ la contraction extr\'emale
associ\'ee \`a $R$. Par l'in\'egalit\'e de Wisniewski, la contraction
est soit de type fibr\'ee soit divisorielle. De plus, le
fibr\'e $L$ est $\phi-$ample et $K_{X}+3L$ supporte la contraction
extr\'emale.\\
\newline
\indent\emph{Etude des contractions de type fibr\'ee.} 
En utilisant \`a nouveau l'in\'egalit\'e
de Wisniewski on v\'erifie que $\text{dim}(Y)\le 3$ et que toute
composante irr\'eductible d'une fibre non triviale est de dimension au
moins $2.$ \\
\newline
\textit{Cas 1 : $\text{dim}(Y)=3.$} Puisque $\phi$ est une contraction
extr\'emale et $\text{dim}(Y)>1$, $\phi$ n'a pas de fibre de
dimension 4. Par 2.(c), le morphisme $\phi$ ne
peut avoir de fibre de dimension $3$ et il en r\'esulte donc que
toutes les fibres de $\phi$ sont de dimension au plus 2, ce qui
entra\^{\i}ne que $Y$ est lisse et que $\phi$ est un fibr\'e projectif
par 1. Remarquons alors
que $L_{|F}\cong\mathcal{O}_{\mathbb{P}^{2}}(1)$ pour toute
fibre $F$ de $\phi$ et consid\'erons la suite exacte :
$$0\longrightarrow\mathcal{T}_{X/Y}\longrightarrow\mathcal{T}_{X} 
\longrightarrow\phi^{*}\mathcal{T}_{Y}\longrightarrow 0$$
La fl\`eche $\mathcal{T}_{X/Y}\longrightarrow L$ obtenue par
composition avec la projection $\mathcal{T}_{X}\longrightarrow L$
\'etant identiquement nulle par le th\'eor\`eme de Grauert et les
r\'esultats ci-dessus, il existe une fl\`eche surjective
$\phi^{*}\mathcal{T}_{Y}\longrightarrow L \longrightarrow 0$ et donc
un morphisme $X\longrightarrow \mathbb{P}_{Y}(\mathcal{T}_{Y})$ au
dessus de $Y$ qui induit un isomorphisme sur chaque fibre. Il en
r\'esulte que ce morphisme est en fait un
isomorphisme, ce qui termine la preuve du th\'eor\`eme dans ce cas.\\
\newline
\textit{Cas 2 : $\text{dim}(Y)=2.$} Par le crit\`ere de
Kobayashi-Ochiai ([9]), une fibre 
g\'en\'erique lisse est une quadrique
$\mathcal{Q}\subset\mathbb{P}^{4}$ de dimension 3. On v\'erifie que 
$L_{|\mathcal{Q}}\cong\mathcal{O}_{\mathcal{Q}}(1)$
et que $H^{0}(\mathcal{Q},\Omega^{1}_{\mathcal{Q}}(1))=0$;
ce cas est \'elimin\'e par le lemme.\\
\newline
\textit{Cas 3 : $\text{dim}(Y)=1.$} Une fibre g\'en\'erique lisse $F$ de
$\phi$ est une vari\'et\'e de Del Pezzo de dimension 4 et $L_{|F}$ est
la polarisation naturelle. En utilisant la classification de T.Fujita
([6], [7], [8]), on v\'erifie que $H^{0}(F,\Omega^{1}_{F}\otimes
L_{|F})=0$ et le lemme  permet de conclure.\\
\newline
\textit{Cas 4 : $\text{dim}(Y)=0.$} Dans ce cas $X$ est de
Fano et, puisque le nombre de Picard de $X$ est 1, on a 
$b_{2}(X)=1$ et on peut appliquer la proposition 1.\\
\newline
\indent\textit{Etude des contractions divisorielles.} 
Notons $E$ le lieu exceptionnel de $\phi$. C'est un
diviseur irr\'eductible. Par l'in\'egalit\'e de Wisniewski, $\phi(E)$
est de dimension 0 ou 1.\\
\newline
\textit{Cas 1 : $\text{dim}(\phi(E))=1.$} Par 2(a), une fibre non
triviale $F$ de $\phi$ est un espace projectif $\mathbb{P}^{3}$ et
$L_{|F}\cong\mathcal{O}_{\mathbb{P}^{3}}(1)$. Ce cas est \`a nouveau
\'eliminer par le lemme.\\
\newline
\textit{Cas 2 : $\text{dim}(\phi(E))=0.$} Dans ce cas, $E$ est soit 
isomorphe \`a $\mathbb{P}^{4}$, soit \`a une
quadrique irr\'eductible de dimension 4, soit \`a  une vari\'et\'e de
Del Pezzo de dimension 4 ([1]). Les deux premiers cas s'\'eliminent
par le lemme. Dans le dernier cas, le fibr\'e normal
$\mathcal{N}_{F|X}$ est $\mathcal{O}_{F}$. Par suite le sch\'ema de
Hilbert $Hilb_{X}$ est lisse au point $F$ et de dimension 1. Puisque
$\phi$ est extr\'emale, les d\'eformations de $F$ doivent \^etre
contract\'ees par $\phi$, ce qui constitue la contradiction
cherch\'ee.\qed\\
\newline 
\textbf{Corollaire}.$-$\textit{Les seules vari\'et\'es de Fano de
dimension 5 admettant une structure de contact sont, \`a isomorphisme
pr\'es, $\mathbb{P}^{5}$,
$\mathbb{P}_{\mathbb{P}^{3}}(\mathcal{T}_{\mathbb{P}^{3}})$ et la
vari\'et\'e de contact homog\`ene de type $G_{2}$.}\\
\newline
\textit{D\'emonstration}.$-$Le corollaire est une cons\'equence
de la conjecture d'Hartshorne-Frankel, d\'emontr\'ee par
S.Mori ([13]).\qed\\
\newline
\indent Pour les vari\'et\'es de dimension de Kodaira $\kappa(X)\ge 0$, nous
avons la\\
\newline 
\textbf{Proposition 2}.$-$\textit{Soit $X$ une vari\'et\'e projective
lisse de dimension $2n+1$. On suppose que $X$ est de dimension de
Kodaira $\kappa(X)\ge 0$. Alors $X$ ne poss\`ede aucune structure de
contact.}\\
\newline
\textit{D\'emonstration}.$-$Raisonnons par l'absurde et supposons que $X$
soit munie d'une structure de contact. Il r\'esulte des  hypoth\`eses,
qu'il existe une vari\'et\'e projective lisse $\overline{X}$ et un
morphisme $\overline{X}\overset{\pi}{\longrightarrow}X$
g\'en\'eriquement fini tel que $h^{0}(\overline{X},\overline{L}^{-1})\ge
1$, o\`u l'on a pos\'e $\overline{L}=\pi^{*}(L)$. Consid\'erons un
ouvert $U\subset X$ non vide au dessus duquel $\pi$ est \'etale et
fini. La structure de contact sur $X$ induit une structure de contact
sur $\pi^{-1}(U)$ associ\'ee au fibr\'e $\overline{L}$. Quitte \`a 
restreindre $\pi^{-1}(U)$, on peut supposer que $\overline{L}$ est
trivialis\'e par une section globale $\overline{\theta}\in
H^{0}(\overline{X},\overline{L}^{-1})\subset
H^{0}(\overline{X},\Omega^{1}_{\overline{X}})$. Sur cet ouvert, la
structure de contact est donn\'ee par la forme de contact
$\overline{\theta}$, ce qui constitue la contradiction cherch\'ee
puisque $d(\overline{\theta})=0$.\qed
\vspace{2cm}\\
\centerline{\textbf{R\'ef\'erences bibliographiques}}
$\ $
\newline
\noindent [1] T.Ando, \emph{On extremal rays of the higher
dimensional varieties}, Invent. Math. 81, 347-357, 1985.\\
\newline
[2] M.Andreatta, J.Wisniewski, \emph{A note on vanishing and
applications}, Duke Math. J. 72, 739-755, 1993.\\
\newline  
[3] M.Andreatta, J.Wisniewski, \emph{A view on
contractions of higher dimensional varieties}, Proc. Sympos. Pure
Math. 62, Part 1, 153-183, 1997.\\
\newline
[4] A.Beauville, \emph{Fano Contact Manifolds and Nilpotent
Orbits}, Comment. Math. Helvet., \`a para\^{\i}tre.\\
\newline
[5] W.M.Boothby, \emph{Homogeneous complex contact manifolds},
Proc. Symp. Proc. Math. 3, 144-154, 1959.\\ 
\newline
[6] T.Fujita, \emph{On the structure of polarized manifods with total
deficeincy one, I}, J. Math. Soc. Japan 32, 709-725, 1980.\\
\newline
[7] T.Fujita, \emph{On the structure of polarized manifods with total
deficeincy one II}, J. Math. Soc. Japan 33, 415-434, 1981.\\
\newline
[8] T.Fujita, \emph{On the structure of polarized manifods with total
deficeincy one III}, J. Math. Soc. Japan 36, 75-89, 1984.\\
\newline
[9] S.Kobayashi, T.Ochiai, \emph{Characterization of complex
projective spaces and hyperquadrics}, J. Math. Kyoto Univ. 13, 31-47, 1973.\\
\newline
[10] Y.Kawamata, K.Matsuda, K.Matsuki, \emph{Introduction to the
minimal model problem}, Adv. Stud. Pure Math. 10, 283-360, 1987.\\
\newline
[11] C.Lebrun, \emph{Fano manifolds, contact structures and
quaternionic geometry}, Int. journ. of Math. 6, 419-437, 1995.\\
\newline
[12] M.Mella, \emph{Existence of good divisors on Mukai manifolds},
alg-geom/9611024, 1996.\\
\newline
[13] S.Mori, \emph{Projective manifolds with ample tangent bundles},
Ann. of Math. 110, 593-606, 1979.\\
\newline
[14] S.Mukai, \emph{Biregular classification of Fano 3-folds and Fano
manifolds of coindex 3}, Proc. Nat. Sci. USA 86, 3000-3002, 1989.\\
\newline
[15] J.Wisniewski, \emph{Length of extremal rays and generalized
adjonction}, Math. Zeit. 200, 409-427, 1989.\\
\newline
[16] J.Wisniewski, \emph{On contractions of extremal rays on Fano
manifolds}, J. reine u. angew Math. 417, 141-157, 1991.\\
\newline
[17] Y-G.Ye, \emph{A note on complex projective threefolds admitting
holomorphic contact structures}, Invent. Math. 121, 421-436, 1995.

\end{document}